\title{FKN, first proof, rewritten}

\author{ Ehud Friedgut and Gil Kalai and Assaf Naor}

\date {}

\documentclass[11pt]{article}
\usepackage{amsmath}
\usepackage{amsthm}
\usepackage{amsfonts}
\numberwithin{equation}{section}

\theoremstyle{plain}
\newtheorem{theo}{Theorem}[section]

\theoremstyle{remark}        

\def\P{{\bf P}}

\def\mu{\P}

\def\le{\leq}
\def\ge{\geq}

\def\hat{\widehat}

\let\CHI=\chi
\def\chi{\raise2.2pt\hbox{$\CHI$}}

\long\def\comment#1{}

\begin{document}
\maketitle

\begin {abstract}
About twenty years ago we wrote a paper, "Boolean Functions whose Fourier Transform is
Concentrated on the First Two Levels", \cite{FKN}. In it we offered several proofs of the statement that 
 Boolean functions $f(x_1,x_2,\dots,x_n)$, whose Fourier coefficients
are concentrated on the lowest two levels are close to a constant function or
to a function of the form $f=x_k$ or $f=1-x_k$. 

Returning to the paper lately, we noticed that the presentation of the first proof is rather cumbersome, and includes several typos. In this note we rewrite that proof, as a service to the public.

\end {abstract}

Here is the main theorem of FKN \cite{FKN}.

\begin {theo}
\label {t:2}
Let $f:\{0,1\}^n \rightarrow \{0,1\}$, 
$\|f\|_2^2 = p$. If $\sum_{|S|>1} \hat f^2(S) = \delta$
then either $p < K' \delta$ or $p > 1- K'\delta$ or
$\|f(x_1, x_2, \dots , x_n)-x_i\|_2^2 \le K \delta$
for some $i$ or
$\|f(x_1, x_2, \dots , x_n)-(1-x_i)\|_2^2 \le K \delta$
for some $i$. Here, $K'$ and $K$ are absolute constants.
\end {theo}

{\bf Proof:} First observe that we may assume that $p=1/2$: replace $f$ by a function $g : \{0,1\}^{n+1} \rightarrow \{0,1\}$ defined by 
$$
g(x_1,\ldots,x_n,x_{n+1}) = f(x_1,\ldots,x_n)\cdot\frac{1+\chi_{n+1}}{2}+(1-f(1-x_1,\ldots,1-x_n))\cdot\frac{1-\chi_{n+1}}{2}
$$
and note that $|g|_2^2=1/2$, and $\sum_{|S|>1} \hat f^2(S) =\sum_{|S|>1} \hat g^2(S)$. Then prove that $g$ is close to a dictator (or anti-dictator). If the influential coordinate for $g$ is $n+1$ then $f$ is necessarily close to a constant, and if the influential coordinate for $g$ is some other $i$, then $f$ too is necessarily close to a dictator or anti-dictator of the $i$th variable.

Let 
$$
f^{=1}: = \sum \hat{f}(i)\chi_i.
$$
Note that $f^{=1}$ is an odd function, i.e. $f^{=1}(S)=-f^{=1}(S^c)$, and that
$$
|f^{=1}|_2^2 = \sum \hat{f}(i)^2= 1/4 - \delta.
$$

Note also that
$$
(1/4-\delta) = \langle f^{=1}, f^{=1}\rangle = \sum \hat{f}(i)^2 = \langle f, f^{=1}\rangle
$$
$$
 =  \sum_{S \subset [n]} 2^{-n}
 f(S) f^{=1}(S) = 
  \sum_{S \subset [n]} 2^{-n}
 f(S) \frac{ f^{=1}(S) - f^{=1}(S^c)}{2} = 
$$
$$
  \sum_{S \subset [n]} 2^{-n}
 f^{=1}(S) \frac{ f(S) - f(S^c)}{2} \le  \sum_{S \subset [n]} 2^{-n}
 \frac{1}{2} \left| f^{=1}(S)\right| =  \frac{1}{2}|f^{=1}|_1.
$$

 So 
 $$
  \frac{1}{2}-2\delta \le |f^{=1}|_1, |f^{=1}|_2 =\sqrt{1/4 -\delta} \le \frac{1}{2}- \delta.
 $$
 
 Normalizing, let $h = \frac{f^{=1}}{|f^{=1}|_2}$, then
 $$
 |h|_2=1, |h|_1 \ge \frac{1/2 -2\delta}{1/2-\delta} \ge 1- 4\delta 
 $$

So $h = \sum a_i \chi_i$, with $\sum a_i^2=1$, and we wish to prove that one of the $a_i$'s is close to 1 or -1.
This follows immediately from the following result of 
K\"onig, Sch\"utt and Tomczak-Jaegermann, \cite{KST}. 
\begin {theo}
\label {t:kst}
Let $a_1,a_2,\dots,a_n$ be real numbers such that $\sum a_i^2 =1$. 
Let $E$ be the expected value of $|\sum a_i\chi_i|$. 
 Then
\begin {equation}
\label {e:kst}
\left|E-\sqrt{\frac{2}{\pi}} \right| \le
\left(1-\sqrt {\frac{2}{\pi}}\right)
\max_{1 \le i \le n}|a_i|.
\end {equation}
\end {theo}
In our case this yields 
$$
\max\{|a_i|\} \ge \frac{1-\sqrt {\frac{2}{\pi}}-4\delta}{1-\sqrt {\frac{2}{\pi}}} > 1-20\delta.
$$

\end {document}